\documentclass[10pt,a4paper]{amsart}
\usepackage[margin=1in]{geometry}

\usepackage{amsmath,amssymb,amsthm}
\usepackage{tikz-cd}
\usepackage{multicol}
\usepackage{aliascnt}
\usepackage{hyperref}
\usepackage[nameinlink,noabbrev,capitalise]{cleveref}
\usepackage[edge length=1cm]{dynkin-diagrams}
\usepackage{comment}
\usepackage{quiver}

\newtheorem{theorem}{Theorem}[section]
\newaliascnt{proposition}{theorem}
\newtheorem{proposition}[proposition]{Proposition}
\aliascntresetthe{proposition}
\newaliascnt{corollary}{theorem}
\newtheorem{corollary}[corollary]{Corollary}
\aliascntresetthe{corollary}
\newaliascnt{lemma}{theorem}
\newtheorem{lemma}[lemma]{Lemma}
\aliascntresetthe{lemma}
\newaliascnt{claim}{theorem}

\aliascntresetthe{claim}
\newaliascnt{fact}{theorem}

\aliascntresetthe{fact}

\newaliascnt{conjecture}{theorem}

\aliascntresetthe{conjecture}
  \theoremstyle{definition}
\newaliascnt{definition}{theorem}
\newtheorem{definition}[definition]{Definition}
\aliascntresetthe{definition}
\newaliascnt{example}{theorem}

\aliascntresetthe{example}
\newaliascnt{remark}{theorem}

\aliascntresetthe{remark}
\newaliascnt{question}{theorem}

\aliascntresetthe{question}
\newaliascnt{notation}{theorem}
\newtheorem{notation}[notation]{Notation}
\aliascntresetthe{notation}
\newaliascnt{observation}{theorem}

\aliascntresetthe{observation}

\crefname{theorem}{theorem}{theorems}
\Crefname{theorem}{Theorem}{Theorems}
\crefname{proposition}{proposition}{propositions}
\Crefname{proposition}{Proposition}{Propositions}
\crefname{corollary}{corollary}{corollaries}
\Crefname{corollary}{Corollary}{Corollaries}
\crefname{lemma}{lemma}{lemmas}
\Crefname{lemma}{Lemma}{Lemmas}
\crefname{claim}{claim}{claims}
\Crefname{claim}{Claim}{Claims}
\crefname{fact}{fact}{facts}
\Crefname{fact}{Fact}{Facts}
\crefname{step}{step}{steps}
\Crefname{step}{Step}{Steps}
\crefname{conjecture}{conjecture}{conjectures}
\Crefname{conjecture}{Conjecture}{Conjectures}
\crefname{definition}{definition}{definitions}
\Crefname{definition}{Definition}{Definitions}
\crefname{example}{example}{examples}
\Crefname{example}{Example}{Examples}
\crefname{remark}{remark}{remarks}
\Crefname{remark}{Remark}{Remarks}
\crefname{question}{question}{questions}
\Crefname{question}{Question}{Questions}
\crefname{notation}{notation}{notations}
\Crefname{notation}{Notation}{Notations}
\crefname{observation}{observation}{observations}
\Crefname{observation}{Observation}{Observations}

\crefformat{equation}{(#2#1#3)}
\Crefformat{equation}{(#2#1#3)}
\crefrangeformat{equation}{(#3#1#4)--(#5#2#6)}
\Crefrangeformat{equation}{(#3#1#4)--(#5#2#6)}
\crefmultiformat{equation}{(#2#1#3)}{ and~(#2#1#3)}{, (#2#1#3)}{ and~(#2#1#3)}
\Crefmultiformat{equation}{(#2#1#3)}{ and~(#2#1#3)}{, (#2#1#3)}{ and~(#2#1#3)}

\newcommand{\BB}{\operatorname{BB}}
\newcommand{\st}{\operatorname{st}}

\title[The Euler characteristic of spherical Bestvina--Brady kernels]{The Euler characteristic of spherical\\ Bestvina--Brady kernels}

\author{Eduardo Flores Mart\'inez}
\address{School of Mathematical Sciences, University of Southampton, Highfield, Southampton, SO17~1BJ, United Kingdom}
\email{\texttt{E.Flores-Martinez@soton.ac.uk}}

\date{October 13, 2026}

\begin{document}
\begin{abstract}
We compute the Euler characteristic of Bestvina--Brady kernels $\BB(\Gamma)$ for all Artin groups of spherical type. Our method involves discrete Morse theory to describe the cellular structure of a classifying space for the Bestvina--Brady kernel, through the action of $\BB(\Gamma)$ on non-integer level-sets of the universal cover of the Salvetti complex. We extend the method of computation of $\chi (\BB (\Gamma))$ to Artin groups that satisfy the $K(\pi,1)$-conjecture and whose nerves are rationally acyclic, and we provide an example by realising the projective plane as the nerve of a Coxeter group.
\end{abstract}
\maketitle

\section{Introduction}

Kernels of maps from Artin groups to $\mathbb{Z}$ have served as examples of groups with prescribed finiteness properties. For instance, finitely generated groups that are infinitely related, and groups of type $\textup{FP}_n$ that are not of type $\textup{FP}_{n+1}$, that is to say, groups $G$ that admit a length $n$, but not length $n+1$, projective resolution of the trivial $\mathbb{Z}G$-module $\mathbb{Z}$ by finitely generated $\mathbb{Z}G$-modules. These examples date back to the work of Stallings and Bieri in the 1960's and 1970's. Bestvina and Brady generalised the setting in which the groups studied by Stallings and Bieri arose by using discrete Morse theory \cite{BB}. Given a right-angled Artin group $A$, its Bestvina--Brady group is defined as the kernel of the epimorphism $A\to\mathbb{Z}$ that sends every standard generator to $1$. The finiteness properties of Bestvina--Brady groups are determined by the homological connectivity of the defining flag complex, with finite presentability being equivalent to simple connectivity.\\ 
The analogous epimorphism is defined for a general Artin group $A$ and we refer to its kernel as a Bestvina--Brady kernel. Such groups have previously been studied, for example, in \cite{ESCARTINMARTINEZ,BLASCOGARCIA2022106984}, giving an analogous homological characterisation of the finiteness properties of the Bestvina--Brady kernel as in the classical case, in which $A$ is a right-angled Artin group.\\ For any Artin group $A$ there is an associated Coxeter group, defined by treating the generators of the Artin group as reflections, i.e. adding the relations $v^2=1$ to the presentation of $A$. The main theme of this paper is to compute the Euler characteristic of Bestvina--Brady kernels, when the ambient Artin group is of spherical type, that is, when its associated Coxeter group is finite.

Our approach uses the Morse function $f:X\to\mathbb{R}$ on the universal cover $X$ of the Salvetti complex induced by the homomorphism to $\mathbb{Z}$. The Bestvina--Brady kernel $H$ acts freely, cellularly and cocompactly on every level-set $L_t=f^{-1}(t)$, inducing a finite cell structure on the quotient $L_t/H$, whose cells correspond to height slices of Coxeter cells. For non-integer $t$, the number of slices contributed by a spherical standard parabolic subgroup equals the length of the longest element of its corresponding Coxeter group.\\ When the nerve is contractible and the Artin group satisfies the $K(\pi,1)$-conjecture, a classical application of the discrete Morse lemma shows $L_t$ as a model for an $EH$. This in particular allows one to compute the Euler characteristic of $H$. In the case of spherical Artin groups, we prove that the Euler characteristic vanishes in the reducible case and we obtain explicit values for every defining irreducible Coxeter diagram; this is contained in \Cref{cor:char-irr-type}. Further, wedges of Coxeter diagrams realise every integer as the Euler characteristic of a Bestvina--Brady kernel with contractible nerve.\\ Finally, we extend the Euler characteristic method of computation to Artin groups satisfying the $K(\pi,1)$-conjecture whose nerves are rationally acyclic. In this setting the level-set quotient exhibits the rational homology, and hence the Euler characteristic of $H$. This generalises the corresponding computation for Bestvina--Brady groups associated to rationally acyclic flag complexes \cite{LearySaadetoglu2011}.
We now state our main result, which follows from \Cref{thm:main} and the basic computations provided by the special cases contained in \Cref{cor:AnBnDnI2}.

\begin{theorem}\label{cor:char-irr-type}
Let $A_\Gamma$ be a spherical Artin group with Coxeter diagram $\Gamma$, and let $\BB(\Gamma)$ be the kernel of the map $A_\Gamma\to\mathbb{Z}$ that sends all generators of $A_\Gamma$ to $1$. If $A_\Gamma$ is not irreducible, then $\chi (\BB(\Gamma))=0$. For the irreducible types, we have the following.
\vspace{-.1cm}

\noindent
\begin{minipage}[t]{.40\linewidth}
\vspace{0pt}
\begin{itemize}
\item $\chi (\BB(A_n))=(-1)^{n-1}$ 
\vspace{.15cm}
\item $\chi (\BB(B_n))=(-1)^{n-1}n$
\vspace{.15cm}
\item $\chi (\BB(D_n))=(-1)^{n-1}(n-2)$
\vspace{.15cm} 
\item $\chi (\BB(I_2^{(n)}))=2-n$
\end{itemize}
\end{minipage}\hfill
\begin{minipage}[t]{.27\linewidth}
\vspace{0pt}
\begin{itemize}
\item $\chi (\BB(H_3))=8$
\vspace{.15cm}
\item $\chi (\BB(H_4))=-42$
\vspace{.15cm}
\item $\chi (\BB(F_4))=-10$
\end{itemize}
\end{minipage}\hfill
\begin{minipage}[t]{.27\linewidth}
\vspace{0pt}
\begin{itemize}
\item $\chi (\BB(E_6))=-7$
\vspace{.15cm}
\item $\chi (\BB(E_7))=16$
\vspace{.15cm}
\item $\chi (\BB(E_8))=-44$
\end{itemize}
\end{minipage}

\end{theorem}

We also note that every integer can be realised as $\chi (\BB (\Gamma))$ for some labelled graph $\Gamma$ with contractible nerve by \Cref{cor:char-wedge} as follows. Let $\Gamma =D_5\vee I_2^{(2)}$ be the labelled graph defined by adjoining the labelled completions by a single vertex and note that $\chi (\BB (\Gamma^{\vee (n-1)}))=n$, which realises all positive integers, while the Coxeter diagrams $A_3$ and $I_2^{(n)}$ realise the rest of the integers for $n\ge 2$.

\subsection{Outline}
In \Cref{sec:2} we state the basic definitions, conventions, constructions and results that will be central to this article, such as the usual naming conventions for the spherical irreducible Coxeter diagrams, a brief description of the Salvetti complex and the definition of Coxeter cells that we will employ.

\Cref{sec:3} deals with the nerve of a Coxeter diagram $\Gamma$ (which we also associate to its Artin group), and the different quotient spaces $Q_t$ obtained from the level-sets of the usual Morse function, defined on the universal cover of the Salvetti complex. We describe the number of cells of $Q_t$ by slicing Coxeter cells at different heights, which allows us to derive \Cref{thm:main} and \Cref{cor:char-irr-type}. After this we prove a vanishing result for the Euler characteristic whenever the underlying Coxeter diagram is disconnected. This section finishes with the computation of the characteristic when a wedge of Coxeter diagrams is considered, the proof of \Cref{thm:main} and computations of its basic cases, corresponding to the Coxeter diagrams $I_2^{(n)}$ and $A_n$.

Lastly, \Cref{sec:4} extends the method of computation of the Euler characteristic of Bestvina--Brady kernels, to the case of rationally acyclic nerves, in contrast to the spherical case treated in \Cref{thm:main,cor:char-irr-type} for which the nerve is contractible. For this we also require the $K(\pi,1)$-conjecture to hold so that the universal cover of the Salvetti complex is contractible. We finish with a worked example on the computation of $\chi (\BB (\Gamma))$, where the nerve of the Coxeter diagram $\Gamma$ is a triangulation of the projective plane.

\section{Basic definitions}\label{sec:2}

\begin{definition}
    Let $\Gamma$ be a finite simple graph with vertices $v_i$ and edges $e_{ij}$ labelled by integers $m_{ij}\in\mathbb{Z}_{\ge 2}$ for some pairs $v_i,v_j$ of vertices. We call $\Gamma$ a labelled graph. The Artin group $A_\Gamma$ has the following defining presentation.
	\[A_\Gamma=\langle v_i\in\Gamma^{(0)}\ \mid\underbrace{v_i v_j v_i\cdots}_{m_{ij}}=\underbrace{v_j v_i v_j\cdots}_{m_{ij}}\text{ whenever }e_{ij}\in\Gamma^{(1)}\rangle .\]
\end{definition}

The addition of the relations $v_i^2$ for all of the vertices of $\Gamma$ yields the defining presentation for the associated Coxeter group $W_\Gamma$.

Whenever $W_\Gamma$ is finite, the Artin group $A_\Gamma$ is called spherical; in this case, we also call $\Gamma$ spherical. A spherical Artin group $A_\Gamma$ decomposes as a product corresponding to a decomposition of $\Gamma$ in which all connecting edges have label $2$. If the group cannot be decomposed non-trivially, it is called irreducible.

\begin{definition}
	A Coxeter diagram $\Gamma$ is a labelled graph with labels $\ge 3$ in which the labels on edges labelled $3$ are omitted from its visual representation. The labelled completion $\overline{\Gamma}$ is defined as the completion of $\Gamma$ to a complete labelled graph in which all new edges have label $2$.
\end{definition}
 
For a labelled graph $\Gamma$, the notation $A_\Gamma$ refers to the Artin group presentation above, with relations only for its edges. When $\Gamma$ is specified to be a Coxeter diagram, we instead use $A_\Gamma$ as shorthand for $A_{\overline{\Gamma}}$, and likewise $W_\Gamma=W_{\overline{\Gamma}}$. Thus missing edges impose no relation for a labelled graph, but represent label $2$ relations for a Coxeter diagram. Subgraphs and subdiagrams used to define standard parabolic subgroups are understood to be induced on their vertex sets.

\begin{definition}
	Let $\Gamma$ be a labelled graph and let $\phi_\Gamma :A_\Gamma\to\mathbb{Z}$ be the epimorphism that sends all of the generators of $A_\Gamma$ to $1\in\mathbb{Z}$. The Bestvina--Brady kernel $\BB(\Gamma)$ is defined as $\ker (\phi_\Gamma)$.
\end{definition}

\newpage

Next we recall the Coxeter diagrams of irreducible type, where $A_n$ and $B_n$ are defined for all $n\ge 1$ and $D_n$ is defined for all $n\ge 3$, with $D_3=A_3$ and $B_1=A_1$.

\begin{figure}[h]
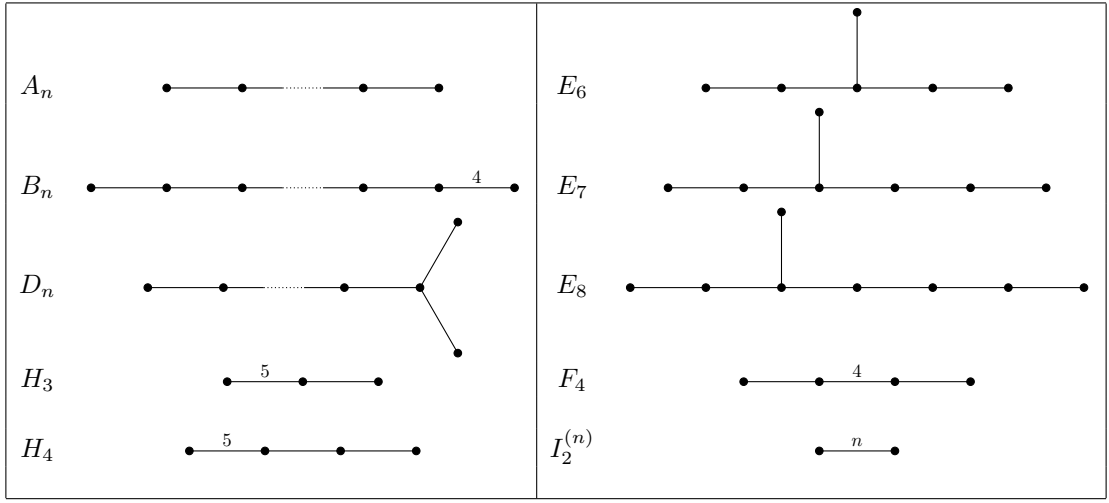

	\centering
	\begin{tabular}{|c c | c c|}
	\hline
	$A_n$ & \dynkin [Coxeter]A{} & $E_6$ & \dynkin [Coxeter]E6\\
	$B_n$ & \dynkin [Coxeter]B{} & $E_7$ & \dynkin [Coxeter]E7 \\
	$D_n$ & \dynkin [Coxeter]D{} & $E_8$ & \dynkin [Coxeter]E8 \\
	$H_3$ & \dynkin [Coxeter]H3 & $F_4$ & \dynkin [Coxeter]F4 \\
	&&&\\
	$H_4$ & \dynkin [Coxeter]H4 & $I_2^{(n)}$ & \dynkin [Coxeter,gonality=n]I{} \\
	&&&\\
	\hline
	\end{tabular}
	\caption{Diagrams of irreducible type.}
	\label{fig:Coxeter-diagrams}
	
\end{figure}

By a classical result of van der Lek \cite{Lek1983TheHT}, for every complete subgraph $S\subseteq\Gamma$, the subgroup $\langle V(S)\rangle$ is isomorphic to $A_S$ and is called a standard parabolic subgroup.

The following description of the Salvetti complex for an Artin group $A_\Gamma$ can be found in \cite{Charney1996FiniteK}, it is built from all of its Coxeter cells as described by Davis \cite{davis2008geometry}. We briefly recall the construction. A finite Coxeter group on $n$ generators acts by orthogonal transformations of $\mathbb{R}^n$, with the generators acting as reflections in codimension-one hyperplanes whose positive half-spaces intersect in a simplicial cone. There exists a point $p$ in the interior of this cone at distance $1$ from every hyperplane. The Coxeter cell associated to a spherical standard parabolic subgroup $A_S$ is defined as the convex hull of the $W_S$-orbit of $p$, and its $1$-skeleton is the Cayley graph of $W_S$. The Salvetti complex is obtained from the wedge of $n$ circles labelled by the generators of $A_\Gamma$ by attaching the Coxeter cells according to the labels on their boundaries.

\begin{notation}\label{notation:longest-length}
Let $\Gamma$ be a labelled graph for which $W_\Gamma$ is finite. We denote by $\ell_\Gamma$ the length of the longest element of $W_\Gamma$, with respect to its standard Coxeter generators. When $\Gamma$ is a Coxeter diagram, we use its labelled completion and set $\ell_\Gamma:=\ell_{\overline{\Gamma}}$.
\end{notation}

\begin{proposition}[\cite{Humphreys_1990}]
The length of the longest element $\ell_\Gamma$ for each spherical irreducible type $\Gamma$, is as listed below.
\vspace{-.15cm}
\begin{multicols}{3}
\begin{itemize}
\item $\ell_{A_n}=n(n+1)/2$
\vspace{.15cm}
\item $\ell_{B_n}=n^2$
\vspace{.15cm}
\item $\ell_{D_n}=n(n-1)$
\vspace{.15cm}
\item $\ell_{I_2^{(n)}}=n$
\end{itemize}

\columnbreak

\begin{itemize}
\item $\ell_{H_3}=15$
\vspace{.15cm}
\item $\ell_{H_4}=60$
\vspace{.15cm}
\item $\ell_{F_4}=24$
\end{itemize}

\columnbreak

\begin{itemize}
\item $\ell_{E_6}=36$
\vspace{.15cm}
\item $\ell_{E_7}=63$
\vspace{.15cm}
\item $\ell_{E_8}=120$
\end{itemize}
\end{multicols}
\end{proposition}

We will be interested in counting the number of heights of a given Coxeter cell. This number equals $\ell_\Gamma$, the length of the longest element of the Coxeter group in the non-integer-height case. For integer heights, it equals $\ell_\Gamma-1$.

\section{Quotient spaces for Bestvina--Brady kernels}\label{sec:3}

We denote by $X$ the universal cover of the Salvetti complex of $A_\Gamma$ and by $X_t$ a level-set of the discrete Morse function $X\to\mathbb{R}$. Such function is defined linearly on each cell of $X$, from its values on $X^{(0)}$ after identifying the vertices of $X$ with the elements of $A_\Gamma$, through the Cayley graph embedding of $A_\Gamma$ in $X$. Note that $X$ is an affine cell complex in the sense of \cite{BB}, where the polyhedral cells associated to the cells of $X$ correspond to the Coxeter cells of $A_\Gamma$. We also let \[Q_t(\Gamma):=X_t/\BB (\Gamma).\] Whenever a quantity depends on a Coxeter diagram $\Gamma$, it is understood to depend on its labelled completion. The topology of $X$ and its local structure at the vertices determine the properties of the level-sets. We recall the relevant definitions of the ascending and descending links below, together with the definition of the nerve.

\begin{definition}\label{def:nerve}
Let $\Gamma$ be a labelled graph. The nerve $N_\Gamma$ is the simplicial complex with vertex set $\Gamma^{(0)}$, whose $m$ simplices are spanned by  subsets $S\subseteq\Gamma^{(0)}$ of cardinality $m+1$ such that the Coxeter group of the induced labelled subgraph on $S$ is finite.
\end{definition}

Thus the faces of the nerve are in correspondence with the spherical standard parabolic subgroups of the Artin group defined by $\Gamma$.

\begin{definition}
	The ascending (respectively, descending) link at $v\in Y^{(0)}$ of a Morse function $f:Y\to\mathbb{R}$ is the full subcomplex of the link of $v$ spanned by the vertices $w$ such that the edge connecting $w$ and $v$ attains its minimum (respectively, maximum) at $v$ under $f$.
\end{definition}

In our case, the ascending and descending links of $X$ are isomorphic to $N_\Gamma$ since $X^{(1)}$ is isomorphic to the Cayley graph of $A_\Gamma$. Thus the $1$-skeleton of each ascending or descending link is isomorphic to $\Gamma$, while its higher-dimensional cells correspond to the spherical subsets of the vertex set of $\Gamma$. This realises $N_\Gamma$ as both the ascending and descending link at every vertex.

Coxeter cells are defined as independent polyhedra. In the following lemma, we identify them with their quotient images in the Salvetti complex and with their lifts to the universal cover.

\begin{lemma}\label{lem:cell-bij}
Let $A_\Gamma$ be an Artin group and $[t]\in (\mathbb{R}-\mathbb{Z})/\mathbb{Z}$. The $m$-cells of $Q_t(\Gamma)$ are in bijection with the $m$-dimensional height slices of the $(m+1)$-dimensional Coxeter cells of the spherical standard parabolic subgroups of $A_\Gamma$.
\end{lemma}
\begin{proof}
 Since all vertices of $X$ have integer height, every $m$-cell of $X_t$ equals the intersection of $X_t$ with a unique $(m+1)$-dimensional Coxeter cell lift. Fix a lift $C$ of the Coxeter cell associated to a spherical standard parabolic subgroup $A_S$, such that it has minimum height equal to $0$. Its maximum height is equal to $\ell_S$, and translating by $g\in A_\Gamma$ shifts all heights by $\phi_\Gamma(g)$. If two translates $gC$ and $hC$ have the same minimum height, then $\phi_\Gamma(g)=\phi_\Gamma(h)$, so $hg^{-1}\in\BB(\Gamma)$ maps one translate to the other while preserving the height slices and their corresponding heights, since $\BB (\Gamma)$ acts cellularly on the level-sets. Thus, each $\BB (\Gamma)$-orbit of an $m$-cell of $X_t$ corresponds to a height slice of its corresponding Coxeter cell lift.
\end{proof}

This lemma is valid for all $[t]\in\mathbb{R}/\mathbb{Z}$, with the distinction that, at integer height all of the $0$-cells are identified so that $Q_t(\Gamma)^{(0)}$ is a single point, and the cell decomposition has one fewer cell for each spherical standard parabolic subgroup. The next figure shows the integer and non-integer height slices of the Coxeter cell, corresponding to the Artin group of type $A_2$, and the corresponding quotient spaces.

\begin{figure}[h]
	\centering
	\includegraphics[width=0.6\textwidth]{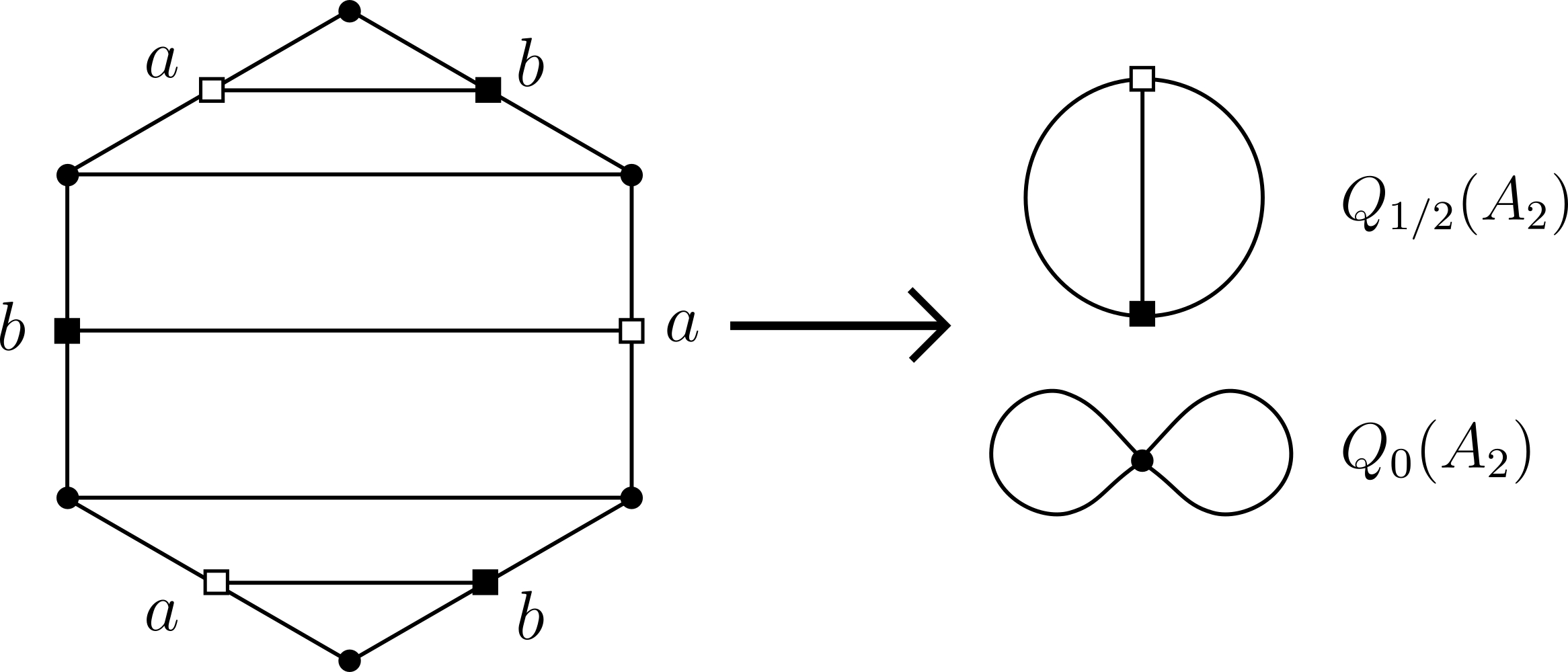}
	\caption{Coxeter cell of type $A_2$ and quotients of level-sets.}
\end{figure}

\begin{corollary}\label{cor:cells-poincare-degree}
Let $\Lambda\subseteq\Gamma$ be a complete subgraph and let $c_m$ ( $\widetilde{c}_m$) be the rank of the $m$-th (reduced) group of the cellular chain complex of $Q_t(\Gamma)$. If $t\in\mathbb{Z}$ then $\widetilde{c}_m=\sum\{\ell_{\Lambda^{(1)}} -1\mid\Lambda\text{ is an }m\text{-face of }N_\Gamma\}$ and if $t\in\mathbb{R}\setminus\mathbb{Z}$ then $c_m=\sum\{\ell_{\Lambda^{(1)}} \mid\Lambda\text{ is an }m\text{-face of }N_\Gamma\}$.
\end{corollary}

\begin{lemma}\label{lemma:spherical-satisfies-Kpi1}
Let $A_\Gamma$ be an Artin group that satisfies the $K(\pi,1)$-conjecture and whose nerve is contractible. Then $Q_t(\Gamma)$ is a $K(\BB (\Gamma),1)$ for all $t\in\mathbb{R}$.
\end{lemma}
\begin{proof}
Since the ascending and descending links of the Morse function on $X$ are isomorphic to $N_\Gamma$, they are contractible. Given that $A_\Gamma$ satisfies the $K(\pi,1)$-conjecture, $X$ is contractible; the Morse lemma \cite{BB} then implies that $X_t$ is also contractible. Since $\BB (\Gamma)$ acts freely on $X_t$, the quotient $Q_t(\Gamma)$ is a classifying space for $\BB (\Gamma)$.
\end{proof}

The following lemma simplifies the formula of $\chi (\BB (\Gamma ))$ in \Cref{thm:main} whenever $A_\Gamma$ is a reducible spherical Artin group.

\begin{lemma}\label{lem:prod-art-semidirect}
    Let $\Gamma'$ and $\Gamma''$ be Coxeter diagrams. Then
    \[\operatorname{BB}(\Gamma'\amalg \Gamma'')\cong (\operatorname{BB}(\Gamma')\times \operatorname{BB}(\Gamma''))\rtimes\mathbb{Z},\]
	in particular if both Bestvina--Brady kernels are of type $F$ then $\chi (\BB (\Gamma'\amalg\Gamma''))$ vanishes. 
\end{lemma}
\begin{proof}
    Let $A'=A_{\Gamma'}$ and $A''=A_{\Gamma''}$. For each Artin group $A_\Lambda\in\{A',A''\}$, consider the short exact sequence
    \[1\to \BB(\Lambda)\to A_\Lambda\xrightarrow{f}\mathbb{Z}\to 1.\]
    Taking the product of the two short exact sequences corresponding to $A'$ and $A''$ and composing with the coordinate sum $\Sigma :\mathbb{Z}\times\mathbb{Z}\to\mathbb{Z}$ gives $\operatorname{BB}(\Gamma'\amalg \Gamma'')=\ker(\Sigma\circ(f'\times f''))$. Projecting $\mathbb{Z}\times\mathbb{Z}$ onto the first coordinate and precomposing with $\left.(f'\times f'')\right|_{\operatorname{BB}(\Gamma'\amalg \Gamma'')}$ gives the surjection $h:\operatorname{BB}(\Gamma'\amalg \Gamma'')\to\mathbb{Z}$ defined by $(a',a'')\mapsto f'(a')$. Since the elements $(a',a'')\in \operatorname{BB}(\Gamma'\amalg \Gamma'')$ satisfy $f'(a')=-f''(a'')$, we have $\ker h=\operatorname{BB}(\Gamma')\times \operatorname{BB}(\Gamma'')$. We therefore obtain the top short exact sequence in the following diagram.
    \[\begin{tikzcd}
	&& {\operatorname{BB}(\Gamma'\amalg \Gamma'')} &&& {\mathbb{Z}} & \\
	1 & {\operatorname{BB}(\Gamma')\times \operatorname{BB}(\Gamma'')} && {A'\times A''} && {\mathbb{Z}\times\mathbb{Z}} & 1 \\
	&&&&& {\mathbb{Z}}
	\arrow[two heads, from=1-3, to=1-6]
	\arrow[from=1-3, to=2-4]
	\arrow[shift left=2, from=1-6, to=2-7]
	\arrow[from=2-1, to=2-2]
	\arrow[from=2-2, to=1-3]
	\arrow[from=2-2, to=2-4]
	\arrow["{f'\times f''}", from=2-4, to=2-6]
	\arrow["{\Sigma\circ (f'\times f'')}"', from=2-4, to=3-6]
	\arrow["{pr_1}", from=2-6, to=1-6]
	\arrow[from=2-6, to=2-7]
	\arrow["\Sigma", from=2-6, to=3-6]
	\arrow[from=3-6, to=2-7]
    \end{tikzcd}\]
When both Bestvina--Brady kernels are of type $F$, we have that $\BB (\Gamma'\amalg \Gamma'')$ can be realised as the fundamental group of a mapping torus with a finite number of cells, built from the finite classifying spaces for $\BB (\Gamma')$ and $\BB (\Gamma'')$. Thus the Euler characteristic of the Bestvina--Brady kernel $\BB (\Gamma'\amalg\Gamma'')$ vanishes since the characteristic telescopically vanishes for a mapping torus.
\end{proof}
The following proposition is a result of Droms, and an alternative proof can be found in \cite{HERMILLER1999141}. We give a proof within the context of classifying spaces.

\begin{proposition}[\cite{DROMS1987484}]\label{prop:kernel-splitting}
	Let $\varphi:A*_CB\to\mathbb{Z}$ be a surjective homomorphism such that the restriction $\left.\varphi\right|_C$ is surjective. Then \[\ker\varphi\cong\ker\left.\varphi\right|_A*_{\ker\left.\varphi\right|_C}\ker\left.\varphi\right|_B.\]
\end{proposition}
\begin{proof}
Let $H=\ker\varphi$ and $H_\bullet=\ker\left.\varphi\right|_\bullet$. Consider the classifying space $B(A*_CB)$ as the pushout $BA\cup_{BC}BB$, and consider the covers associated to $A$, $B$ and $C$. These covers are classifying spaces of the corresponding subgroups, since covering spaces induce isomorphisms on higher homotopy groups. As the restrictions of $\varphi$ to $A,B$ and $C$ are surjective, the associated covers are connected. Taking the preimages of $BA$, $BC$ and $BB$ under the covering map $BH\to B(A*_CB)$ yields $BH\simeq BH_A\cup_{BH_C}BH_B$, and the result follows from van Kampen's theorem.
\end{proof}
If, instead of the disjoint union in \Cref{lem:prod-art-semidirect}, we consider the wedge of two Coxeter diagrams, we obtain the following corollary, which allows us to realise every integer as the Euler characteristic of a Bestvina--Brady kernel. Given two Coxeter diagrams $\Gamma$ and $\Lambda$, we denote by $\Gamma\vee\Lambda$ the labelled graph obtained by adjoining the labelled completions $\overline{\Gamma}$ and $\overline{\Lambda}$ at a single vertex.

\begin{corollary}\label{cor:char-wedge}
	Let $\Gamma$ and $\Lambda$ be Coxeter diagrams whose Bestvina--Brady kernels are of type $F$, and denote by $\Gamma\vee\Lambda$ their wedge at a single vertex. Then \[\chi (\BB (\Gamma\vee\Lambda ))=\chi (\BB(\Gamma))+\chi(\BB (\Lambda))-1.\]
\end{corollary}
\begin{proof}
	By \Cref{prop:kernel-splitting}, we have $\BB (\Gamma\vee\Lambda)\cong\BB (\Gamma)*\BB (\Lambda )$, and the result follows from the usual Euler characteristic formula for a wedge.
\end{proof}

\begin{theorem}\label{thm:main} Let $A_\Gamma$ be a spherical Artin group on $n$ generators with defining Coxeter diagram $\Gamma$. If $A_\Gamma$ is not irreducible, then $\chi (\BB (\Gamma))=0$. If $A_\Gamma$ is irreducible, then
    \[\chi (\BB (\Gamma))=(-1)^{n-1}\ell_\Gamma+\sum_{k=1}^{n-1}\sum_{\Lambda\in\mathcal{C}_{n-k}}(-1)^{k-1}\chi (\BB (\Lambda ))\]
    where $\mathcal{C}_m:=\{\Lambda\subseteq\Gamma :\Lambda\text{ is connected and }|\Lambda^{(0)}|=m\}$.
\end{theorem}
\begin{proof}
	All spherical Artin groups satisfy the $K(\pi,1)$-conjecture \cite{Deligne1972}, and their nerves are simplices. Thus, by \Cref{lemma:spherical-satisfies-Kpi1}, the Euler characteristic of $\BB (\Gamma)$ equals $\chi (Q_t(\Gamma))$ for $t\in\mathbb{R}\setminus\mathbb{Z}$. By \Cref{cor:cells-poincare-degree}, the space $Q_t(\Gamma)$ has $\ell_\Gamma$ top-dimensional cells of dimension $n-1$. Therefore, by the definition of Euler characteristic,
	\begin{equation}\label{eq:char-minus-topcells}
	\chi (\BB (\Gamma))=(-1)^{n-1}\ell_\Gamma +\chi (Q_t(\Gamma)^{(n-2)}).
	\end{equation}
	Let $V=\Gamma^{(0)}$, and for $J\subseteq V$ let $\Gamma_J$ denote the induced Coxeter subdiagram on $J$. For each $v\in V$, consider the embedded subcomplex $Y_v:=Q_t(\Gamma_{V\setminus\{v\}})\subseteq Q_t(\Gamma)$.
	Each cell of $Q_t(\Gamma)^{(n-2)}$ is indexed by a proper subset of $V$, and is therefore contained in at least one $Y_v$. Hence
	\[
	Q_t(\Gamma)^{(n-2)}=\bigcup_{v\in V}Y_v.
	\]
	If $I\subsetneq V$ is not empty, then a cell belongs to every $Y_v$ with $v\in I$ precisely when its indexing subset is contained in $V\setminus I$. Consequently,
	\[
	\bigcap_{v\in I}Y_v=Q_t(\Gamma_{V\setminus I}).
	\]
	Note that $\bigcap_{v\in V}Y_v$ is empty, since cells are indexed by nonempty subsets of $V$. Thus a nonempty $k$-fold intersection is indexed by a unique set $I\subsetneq V$ of $k$ deleted vertices and corresponds to the induced Coxeter subdiagram on the remaining $n-k$ vertices. Hence inclusion--exclusion and \Cref{lemma:spherical-satisfies-Kpi1} give
	\begin{equation}\label{eq:char-of-n-2-sk}
	\chi(Q_t(\Gamma)^{(n-2)})
	=\sum_{k=1}^{n-1}(-1)^{k-1}
	\sum_{\substack{J\subseteq V\\|J|=n-k}}
	\chi(\BB(\Gamma_J)).
	\end{equation}
	A spherical Coxeter diagram is irreducible if and only if it is connected. Moreover, by \Cref{lem:prod-art-semidirect}, the Euler characteristic of the Artin kernel associated to a reducible spherical diagram vanishes. Therefore, the only contributing Coxeter subdiagrams are those in $\mathcal{C}_{n-k}$. Substituting \Cref{eq:char-of-n-2-sk} into \Cref{eq:char-minus-topcells} yields the result.
\end{proof}

We now apply \Cref{thm:main} and the corresponding lengths of the longest elements \cite{Humphreys_1990} to spherical irreducible types to obtain \Cref{cor:char-irr-type}. All cases other than $A_n$ and $I_2^{(n)}$ follow from these two irreducible types, together with the computation of $\ell_\Gamma$ in \cite{Humphreys_1990} via the number of roots, as shown in the next diagram, where an arrow from a diagram $\Gamma$ to a diagram $\Lambda$ means that the computation of $\chi (\BB (\Lambda))$ depends on the value of $\chi (\BB (\Gamma))$. We briefly compute the Euler characteristic for the families $A_n$, $B_n$ and $D_n$.

\begin{figure}[h]
\[\begin{tikzcd}[ampersand replacement=\&]
	\& {I_2^{(5)}} \& {H_3} \& {H_4} \& {E_6} \& {E_7} \& {E_8} \& \\
	\&\& {D_4} \& {D_5} \& {D_6} \& {D_7} \& {D_8} \& \cdots \\
	{F_4} \& {A_2} \& {A_3} \& {A_4} \& {A_5} \& {A_6} \& {A_7} \& \cdots \\
	\& {I_2^{(4)}} \& {B_3} \& {B_4} \& \cdots
	\arrow[from=1-2, to=1-3]
	\arrow[from=1-3, to=1-4]
	\arrow[from=1-5, to=1-6]
	\arrow[from=1-6, to=1-7]
	\arrow[from=2-3, to=2-4]
	\arrow[from=2-4, to=1-5]
	\arrow[from=2-4, to=2-5]
	\arrow[from=2-5, to=1-5]
	\arrow[from=2-5, to=1-6]
	\arrow[from=2-5, to=2-6]
	\arrow[from=2-6, to=1-6]
	\arrow[from=2-6, to=1-7]
	\arrow[from=2-6, to=2-7]
	\arrow[from=2-7, to=1-7]
	\arrow[from=2-7, to=2-8]
	\arrow[from=3-2, to=1-3]
	\arrow[from=3-2, to=3-3]
	\arrow[from=3-2, to=4-3]
	\arrow[curve={height=-40pt}, from=3-3, to=1-4]
	\arrow[from=3-3, to=2-3]
	\arrow[from=3-3, to=3-4]
	\arrow[from=3-3, to=4-4]
	\arrow[from=3-3, to=4-5]
	\arrow[from=3-4, to=2-4]
	\arrow[from=3-4, to=3-5]
	\arrow[curve={height=-12pt}, from=3-5, to=1-5]
	\arrow[from=3-5, to=2-5]
	\arrow[from=3-5, to=3-6]
	\arrow[curve={height=-12pt}, from=3-6, to=1-6]
	\arrow[from=3-6, to=2-6]
	\arrow[from=3-6, to=3-7]
	\arrow[curve={height=-12pt}, from=3-7, to=1-7]
	\arrow[from=3-7, to=2-7]
	\arrow[from=3-7, to=3-8]
	\arrow[from=4-2, to=4-3]
	\arrow[from=4-3, to=3-1]
	\arrow[from=4-3, to=4-4]
	\arrow[from=4-4, to=4-5]
\end{tikzcd}\]
\caption{Inclusions of Coxeter diagrams.}
\end{figure}

\begin{corollary}\label{cor:AnBnDnI2}
For the spherical Artin groups of types $A_n$, $B_n$, $D_m$, and
$I_2^{(k)}$, the Euler characteristics of the associated
Bestvina--Brady kernels are as follows, where $n\ge 1$, $m\ge 3$ and $k\ge 2$.
\vspace{-.15cm}
\begin{multicols}{2}
\begin{itemize}
\item $\chi(\BB(A_n)) = (-1)^{n-1}$
\vspace{.15cm}
\item $\chi(\BB(B_n)) = (-1)^{n-1}n$
\end{itemize}
\columnbreak
\begin{itemize}
\item $\chi(\BB(D_m)) = (-1)^{m-1}(m-2)$
\vspace{.15cm}
\item $\chi(\BB(I_2^{(k)})) = 2-k$
\end{itemize}
\end{multicols}
\end{corollary}
\begin{proof}
The case $I_2^{(n)}$ is immediate from the fact that $\chi (pt)=1$ and that the corresponding maximal Coxeter cell is a regular $2n$-gon, which implies that $\ell_{I_2^{(n)}}=n$. The Euler characteristic of $\BB(I_2^{(n)})$ can\vspace{-.15cm}\\ also be seen directly to be $2-n$ since for an integer level-set, the space $Q_t(I_2^{(n)})$ is a wedge of $n-1$ circles.
    By strong induction and \Cref{thm:main}, we have the following recursive formulas, from which the corresponding values follow.

	\begin{align*}
	\chi(\BB(A_{n+1}))&=(-1)^n\ell_{A_{n+1}}+\sum_{k=2}^{n+1}(-1)^kk\chi (\BB (A_{n+2-k})),\\
	\chi (\BB (B_{n+1}))&=(-1)^n(n+1)^2+\sum_{k=0}^{n-1}(-1)^k(\chi (\BB (B_{n-k}))+(k+1)\chi (\BB (A_{n-k}))),\\
	\chi (\BB(D_{n+1}))&=(-1)^nn(n+1)+\sum_{k=0}^{n-3}(-1)^k\chi (\BB (D_{n-k}))+\sum_{k=0}^{n-1}(-1)^k((k+2)\chi (\BB (A_{n-k}))).\qedhere
	\end{align*}
\end{proof}

\section{Rationally acyclic nerve}\label{sec:4}

We now extend our method of computation of $\chi (\BB (\Gamma))$ to Artin groups that satisfy the $K(\pi,1)$-conjecture, and whose nerve is rationally acyclic. We give an explicit example which arises from the minimal triangulation of the projective plane.

\begin{lemma}\label{lemma:Q-acyclic-nerve}
	Let $A_\Gamma$ be a finitely generated Artin group that satisfies the $K(\pi,1)$-conjecture and whose nerve is $\mathbb{Q}$-acyclic. Then $H_*(\BB(\Gamma);\mathbb{Q})\cong H_*(Q_t(\Gamma);\mathbb{Q})$ for every $t\in\mathbb{R}$.
\end{lemma}
\begin{proof}
The ascending and descending links of the Morse function induced by $\phi_\Gamma$ are isomorphic to the nerve of $\Gamma$. Since the universal cover of the Salvetti complex is contractible by the $K(\pi,1)$-conjecture, the Morse lemma, as in \cite{BB}, implies that the level-sets $X_t$ are $\mathbb{Q}$-acyclic. The group $\BB(\Gamma)$ acts freely and cellularly on $X_t$, so the cellular $\mathbb{Q}$-chain complex of $X_t$ is a finite free resolution of $\mathbb{Q}$. Since $\mathbb{Q}\otimes_{\mathbb{Q}\BB(\Gamma)}C_*(X_t;\mathbb{Q})\cong C_*(Q_t(\Gamma);\mathbb{Q})$, the result follows.
\end{proof}
This weakened version of \Cref{lemma:spherical-satisfies-Kpi1} implies that the Euler characteristic of $\BB (\Gamma)$ can be computed from the cell count of $Q_t(\Gamma)$. Next, we describe the stellar subdivision of a simplicial complex, which allows us to realise known $\mathbb{Q}$-acyclic spaces as nerves of Coxeter groups.

\begin{definition}
	Let $\Gamma$ be a labelled graph and $\Delta$ a complete labelled subgraph. The stellar subdivision $\st(\Gamma,\Delta)$ is defined as the labelled graph obtained by adding a new vertex to $\Gamma$ and connecting it to the vertices of $\Delta$ with label $2$ edges.
\end{definition}

\begin{lemma}\label{lemma:K-pi-one}
	Let $\Gamma$ be a labelled graph and $\Delta$ a complete subgraph. Then $A_{\st(\Gamma, \Delta)}$ satisfies the $K(\pi,1)$-conjecture if $A_\Gamma$ and $A_\Delta$ satisfy it as well.
\end{lemma}
\begin{proof}
	Let $S(A)$ denote the Salvetti complex associated to the Artin group $A$. Since $\st(\Gamma,\Delta)$ is the adjunction of $\Gamma$ with a star graph whose edges all have label $2$, the space $S(A_{\st(\Gamma, \Delta)})$ is the adjunction space of $S(A_\Gamma)$ and $S(A_\Delta)\times S^1$ through $S(A_\Delta)$. Since a classifying space for an amalgamated product is the pushout of the corresponding classifying spaces, we have $S(A_{\st(\Gamma,\Delta)})=BA_{\st(\Gamma,\Delta)}$.
\end{proof}

We now apply the previous results to compute the Euler characteristic of a Bestvina--Brady kernel that is a subgroup of an Artin group satisfying the $K(\pi,1)$-conjecture and having a non-contractible, rationally acyclic nerve.

\begin{figure}[h]
	\centering
	\includegraphics[scale=.15]{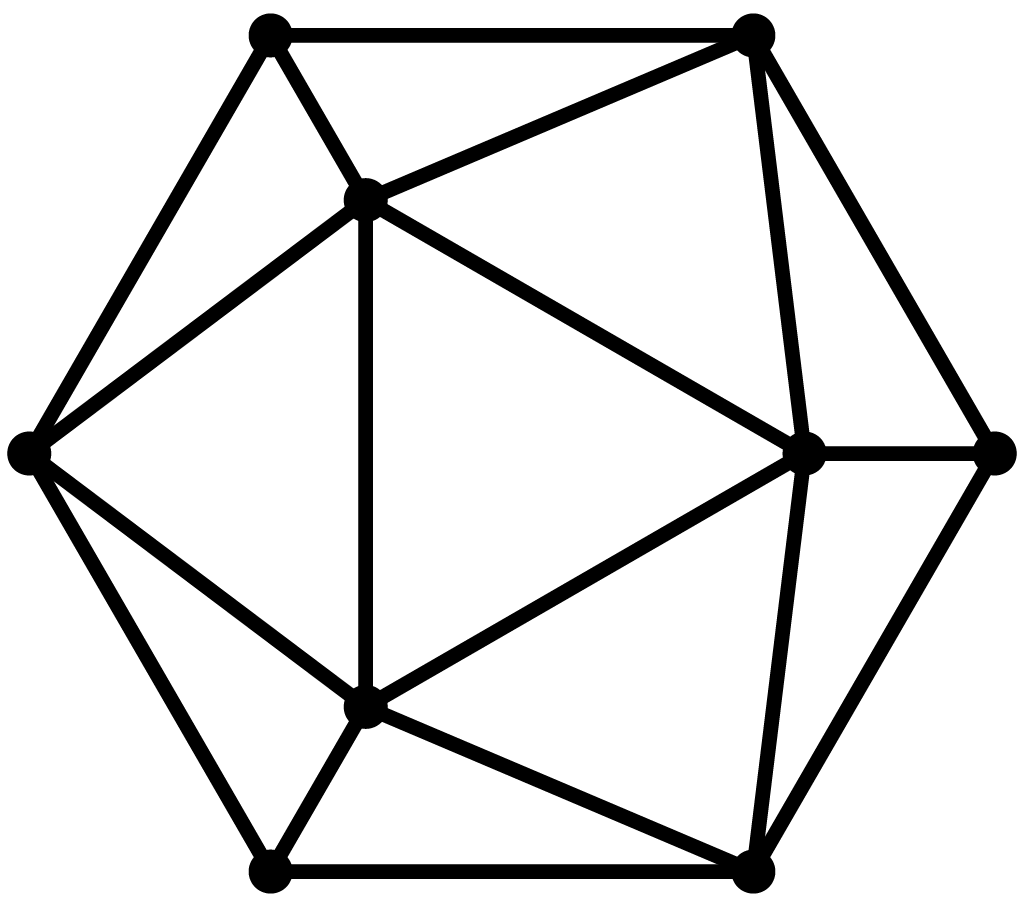}
	\caption{Minimal triangulation of $\mathbb{R}\mathbb{P}^2$.}
	\label{fig:min-tr-RP2}
\end{figure}

Let $T$ be the minimal $6$-vertex triangulation of the projective plane depicted in \Cref{fig:min-tr-RP2}, before the antipodal identification of the boundary edges, and label each edge of $\Gamma:=T^{(1)}$ by $3$. Since $T^{(1)}$ is a labelled graph and large-type Artin groups satisfy the $K(\pi,1)$-conjecture \cite{hendriks1985hyperplane}, the group $A_\Gamma$ satisfies the conjecture. Let $\st (\Gamma)$ be the iterated stellar subdivision of $\Gamma$ over all of the $10$ triangles, each of which consists of the vertices of a $2$-simplex of the $\mathbb{R}\mathbb{P}^2$ triangulation. It follows from \Cref{lemma:K-pi-one} that the Artin group $A_{\st(\Gamma)}$ also satisfies it.

Each of the $30$ triangles in $\st(\Gamma)$ has the labelling $223$ and generates a spherical Artin group with Coxeter diagram $A_2\amalg A_1$; moreover, there are no parabolic subgroups with more than three generators. Thus the nerve $N_{\st(\Gamma)}$ is a triangulation of the projective plane. By \Cref{lemma:Q-acyclic-nerve}, we can compute $\chi (\BB(\st(\Gamma)))$ by counting the cells in $Q_t(\st(\Gamma))$. We have $\ell_{A_2\amalg A_1}=4$ since, for any Coxeter diagram $\Lambda$ the quantity $\ell_\Lambda$ equals the degree of its Poincar\'e polynomial, and the degree associated to a disjoint union of Coxeter diagrams is the sum of their degrees. Since $T$ has $10$ two-dimensional faces, \Cref{cor:cells-poincare-degree} gives, for $t\in\mathbb{Z}$,
\[\sum\{\ell_{\Lambda^{(1)}} -1\mid\Lambda\text{ is a }2\text{-face of }N_{\st(\Gamma)}\}=90\]
and
\[\chi(\BB(\st(\Gamma)))=90+\chi (Q_t(\st (\Gamma))^{(1)})=90+\chi (\vee_{i=1}^{2\cdot 15+30}S^1)=90-59=31.\]

\section*{Acknowledgements}
I am grateful to Ian Leary for many helpful discussions concerning Bestvina--Brady kernels, comments on earlier versions of this manuscript and the creation of \Cref{fig:min-tr-RP2}. I would also like to thank Ashot Minasyan for his comments on this project, and in particular for suggesting the use of the length of the longest element of a Coxeter group in place of the degree of its Poincar\'e polynomial.

\section*{AI Declaration}
OpenAI Prism with GPT-5.2 was employed for language editing, \LaTeX{} formatting, notational consistency checks and revision of the exposition. The author reviewed all AI-assisted suggestions incorporated into the manuscript and takes full responsibility for its contents.

\bibliographystyle{amsalpha}
\bibliography{references}
\end{document}